\documentclass[a4paper,reqno]{amsart}
\usepackage{amssymb, amsmath, amscd}
\def\XIx\langle#1\rangle{h(#1)}
\newtheorem{theorem}{Theorem}[section]
\newtheorem{definition}[theorem]{Definition}
\newtheorem{lemma}[theorem]{Lemma}

\newtheorem{remark}[theorem]{Remark}

\makeatletter
 \@addtoreset{equation}{section}
\makeatother
\begin{document}
\title[Homogeneity of Lorentzian manifolds with recurrent curvature]
{Homogeneity of Lorentzian three-manifolds with recurrent curvature}
\author{E. Garc\'{\i}a-R\'{\i}o, P.  Gilkey  \text{and} S. Nik\v cevi\'c}
\address{EG: Faculty of Mathematics, University of Santiago de Compostela, Spain}
\email{eduardo.garcia.rio@usc.es}
\address{PG: Mathematics Department, \; University of Oregon, \;\;
  Eugene \; OR 97403 \; USA}
\email{gilkey@uoregon.edu}
\address{SN: Mathematical Institute, Sanu, Knez Mihailova 36, p.p. 367,
11001 Belgrade, Serbia}
\email{stanan@mi.sanu.ac.rs}


\begin{abstract}
$k$-Curvature homogeneous three-dimensional Walker metrics are described for $k\leq 2$.
This allows a complete description of locally homogeneous three-dimensional Walker metrics, showing that there exist exactly three isometry classes of such manifolds.
As an application one obtains a complete description of all locally homogeneous Lorentzian manifolds with recurrent curvature.
Moreover, potential functions are constructed in all the locally homogeneous manifolds resulting  in steady gradient Ricci and Cotton solitons.
\end{abstract}
\maketitle

\section{Introduction}\label{sect-1}

The study of Lorentzian three-manifolds admitting a parallel null vector field is central both in geometry and physics. Physically they represent the simplest non-trivial pp-waves and, from a geometrical point of view, they are the underlying structure of many Lorentzian situations without Riemannian counterpart.
A Lorentzian manifold
is said to be \emph{irreducible} if the holonomy group does not leave invariant any non-trivial subspace.
Moreover, the action of the holonomy is said to be \emph{indecomposable} if it leaves invariant only
non-trivial subspaces for which the restriction of the metric degenerates. Then the de Rham-Wu's Theorem shows that any complete and simply connected Lorentzian manifold is a product of indecomposable ones.
Thus Walker three-manifolds constitute the basic material to build many Lorentzian metrics.

When trying to analyze the curvature of a given manifold, one must deal not only with the curvature tensor itself, but also with its covariant derivatives. In the locally symmetric case, the curvature tensor is parallel and hence the study can be reduced to the purely algebraic level. Generalizing this condition, one naturally considers the case when the curvature tensor is {\it recurrent} (i.e., $\nabla R=\omega\otimes R$ for some $1$-form $\omega$) but not parallel ($\nabla R\neq 0$). The class of recurrent Lorentzian manifolds reduces to the study of plane waves and the three-dimensional Walker manifolds \cite{galaev-recurrent, walker-recurrent}.
Homogeneous plane waves were discussed in \cite{BL} and hence one of the purposes of this paper is to give a complete description of all locally homogeneous Walker three-manifolds.

Our main result shows that there exist exactly three isometry classes of such manifolds (cf. Theorem \ref{thm-1.4}), which allows a complete description of all locally homogeneous Lorentzian manifolds with recurrent curvature.
Recall that, in any locally homogeneous pseudo-Riemannian manifold, the curvature tensor and all its covariant derivatives are the same at each point $p$ of the manifold. Generalizing this condition, a manifold $(M,g)$ is said to be \emph{$k$-curvature homogeneous} if for any two points there exists a linear isometry between the corresponding tangent spaces which preserves the curvature tensor and its covariant derivatives up to order $k$. Clearly any locally homogeneous space is curvature homogeneous of any order; conversely, given $m$, there is $k=k(m)$ so that any
$k$-curvature homogeneous manifold is in fact homogeneous \cite{PS,Singer}.
Our approach is based on a careful analysis of the curvature homogeneity of a Walker three-manifold.

Here it is worth emphasizing that Lorentzian manifolds which are $1$-curvature homogeneous play a distinguished role (see, for example \cite{bueken-djoric}), and it is also a purpose of this work to give a complete description of all Walker three-manifolds which are $1$-curvature homogeneous (cf. Theorem \ref{thm-1.3}).
Lorentzian manifolds of dimension $3$ which are $1$-curvature homogeneous are as close as possible to being locally homogeneous and play a basic role in constructing many interesting Lorentzian examples without Riemannian counterpart \cite{Ca, CLGRVAVL}. It is our second purpose in this paper to analyze their geometry, thus showing  that any $1$-curvature homogeneous Walker three-manifold is either a gradient Ricci soliton or a gradient Cotton soliton.

Recall that, generalizing the Einstein condition, a Lorentzian manifold $(M,g)$ is said to be a {\it Ricci soliton} if and only if it is a self-similar solution of the Ricci flow, i.e., the one-parameter family of metrics $g(t)=\sigma(t)\psi^*_t g$ is a solution of the Ricci flow $\frac{\partial}{\partial t}g(t)=-2 Ric_{g(t)}$, for some smooth function $\sigma(t)$ and some one-parameter family of diffeomorphisms $\psi_t$ of $M$, where $Ric_{g(t)}$ denotes the Ricci tensor of $(M, g(t))$.
From a physical viewpoint, Ricci solitons can be interpreted as special solutions of the Einstein field equations, where the stress-energy tensor essentially corresponds to the Lie derivative of the metric.
We analyze in detail the structure of gradient Ricci solitons on Walker three-manifolds in Section \ref{se:rs}, showing that any $1$-curvature homogeneous Lorentzian three-manifold with recurrent curvature is indeed a steady gradient Ricci soliton. While one of the possible $1$-curvature homogeneous Walker three-manifold is indeed  a plane wave and thus a expanding, steady and shrinking Ricci soliton, it is shown (cf. Theorem~\ref{thm-2.7}) that the non-homogeneous family does not support any non-steady Ricci soliton.

The Cotton tensor, $C$, measures the failure of the Schouten tensor to be Codazzi. The existence of self-similar solutions to the Cotton flow
$\frac{\partial}{\partial t}g(t)=-\lambda C_{g(t)}$ also provides a family of three-dimensional metrics which generalize the locally conformally flat manifolds. Locally homogeneous Walker three-manifolds split essentially into two families, one of which is locally conformally flat.
The existence of gradient Cotton solitons is also studied, showing that any locally homogeneous Walker three-manifold is a steady gradient Cotton soliton.  Moreover, a locally homogeneous Walker three-manifold admits a non-steady Cotton soliton if and only if it is locally conformally flat as a consequence of Theorem~\ref{th-1}.

\section{Preliminaries}

\subsection{Walker metrics}

A Lorentzian manifold admitting a parallel null vector field will be refereed in what follows as a {\it Walker manifold}
(following the notation in \cite{walker-metrics}).
A three-dimensional Walker metric admits local adapted coordinates $(x,y,\tilde x)$
where the metric is given by
\begin{equation}\label{eq:walker-metric}
g(\partial_x,\partial_x)=-2f(x,y),\quad g(\partial_x,\partial_{\tilde x})
=g(\partial_y,\partial_y)=1\,.
\end{equation}
We shall work locally. Let $\mathcal{O}$ be
a connected open subset of $\mathbb{R}^2$, let $M:=\mathcal{O}\times\mathbb{R}$,
and let $(x,y,\tilde x)$ on $\mathbb{R}^3$ be coordinates on $M$.
For $f\in C^\infty(\mathcal{O})$, let $\mathcal{M}_f:=(M,g_f)$
where $g=g_f$ is the Lorentz metric on $M$ given by \eqref{eq:walker-metric}.
We note for future reference that the non-zero covariant derivatives are given by:
$$
\nabla_{\partial_x}\partial_x=-f_x\partial_{\tilde x}+f_y\partial_{y}
\quad\text{and}\quad
\nabla_{\partial_x}\partial_y=\nabla_{\partial y}\partial_x=-f_y\partial_{\tilde x}.
$$
\begin{remark}\rm
Since $\nabla_x\partial_x=-f_x\partial_{\tilde x}+f_y\partial_{y}$
involves a $\partial_y$ component which can depend on $y$,
these manifolds are not generalized plane wave manifolds and in fact
can be geodesically incomplete; they can exhibit Ricci blowup -- see
Section 3.2 of \cite{G07} for further details.
\end{remark}

The only (potentially) non-zero components of the $(1,3)$ curvature tensor $\mathcal{R}$  are given by
$$
\mathcal{R}(\partial_x,\partial_y)\partial_y=f_{yy}\partial_{\tilde x}
\quad
\mbox{and}
\quad
\mathcal{R}(\partial_x,\partial_y)\partial_x=-f_{yy}\partial_y.
$$
Thus the $(0,4)$ curvature tensor is determined by
$R(\partial_x,\partial_y,\partial_y,\partial_x)=f_{yy}$.
Similarly, when considering $\nabla R$, the only possible contributions
from the Christoffel symbols arise when plugging $\partial_x$ or $\partial_y$ in
the last entry
$\nabla R(\partial_x,\partial_y,\partial_y,\partial_x;\,\cdot\,)$.
Consequently the (possibly) non-zero entries in $\nabla R$ are given by:
$$
\nabla R(\partial_x,\partial_y,\partial_y,\partial_x;\partial_x)=f_{xyy}
\quad
\mbox{and}
\quad
\nabla R(\partial_x,\partial_y,\partial_y,\partial_x;\partial_y)=f_{yyy}.
$$
\begin{remark}\rm
It follows from previous equations that $\mathcal{M}_f$ is locally symmetric if and only if $f_{yy}$ is a constant. Whenever $f_{yy}=\operatorname{const.}\neq 0$, the resulting manifold is a Cahen-Wallach symmetric space \cite{CLPTV}. In dimension $3$, all Walker metrics have recurrent curvature in a neighborhood of any point of non-zero curvature, i.e., $\nabla R=\omega\otimes R$, for a $1$-form
$$
\omega= (\ln f_{yy})_x \,\, dx +  (\ln f_{yy})_y \,\,dy .
$$
We refer to \cite{walker-metrics} and the references therein for more information on Walker three-manifolds.
\end{remark}
We compute similarly that the (possibly) non-zero entries in  $\nabla^2R$ are:
$$
\begin{array}{l}
\nabla^2R(\partial_x,\partial_y,\partial_y,\partial_x;\partial_x,\partial_x)
=f_{xxyy}-f_yf_{yyy},\\
\noalign{\medskip}
\nabla^2R(\partial_x,\partial_y,\partial_y,\partial_x;\partial_x,\partial_y)
=f_{xyyy},\\
\noalign{\medskip}
\nabla^2R(\partial_x,\partial_y,\partial_y,\partial_x;\partial_y,\partial_x)
=f_{xyyy},\\
\noalign{\medskip}
\nabla^2R(\partial_x,\partial_y,\partial_y,\partial_x;\partial_y,\partial_y)
=f_{yyyy}\,.
\end{array}
$$
Next we introduce the following special classes of metrics \eqref{eq:walker-metric} which play a distinguished role in our analysis.
\begin{definition}\label{def-2.3}\rm
Let $\mathcal{M}_f$ be a three-dimensional Walker manifold given by \eqref{eq:walker-metric}.
\begin{enumerate}
\item Let $\mathcal{N}_b$ be defined by taking $f(x,y)=b^{-2}e^{by}$ for $b\ne0$.
\item Let $\mathcal{P}_c$ be defined by taking $f(x,y)=\frac12y^2\alpha(x)$
where $\alpha_x=c\alpha^{3/2}$ and $\alpha>0$.
\item Let $\mathcal{CW}_\varepsilon$ be the three-dimensional Cahen-Wallach symmetric space defined by taking $f(x,y)=\varepsilon y^2$.
%Note. We may take
%$\alpha(x)=-\tilde c(x-x_0)^{-2}$ for $x>x_0$ if $c<0$ and we may take
%$\alpha(x)=\tilde c(x-x_0)^{-2}$
%if $c>0$ and $x<x_0$ where $\tilde c>0$ is suitably chosen.
\end{enumerate}
\end{definition}

Then one has the following result, which follows from theorems \ref{thm-1.3} and \ref{thm-1.4} together with the  results in \cite{GN05},

 \begin{theorem}\label{thm-1.2}
 \ \begin{enumerate}
\item The manifolds $\mathcal{CW}_\varepsilon$ are locally symmetric.
\item The manifolds $\mathcal{N}_b$ and $\mathcal{P}_c$ are locally
homogeneous.
\item The manifolds $\{\mathcal{CW}_\varepsilon,\mathcal{N}_b,\mathcal{P}_c\}$
 have non-isomorphic $1$-curvature models and represent different
 local isometry types.
\end{enumerate}
\end{theorem}

\begin{remark}\rm
 The manifolds $\mathcal{P}_c$ are plane-waves (as well as the manifolds $\mathcal{CW}_\varepsilon$) since the function $f(x,y)$ is quadratic in $y$; they are not generalized plane wave manifolds as discussed in \cite{GN05} and we regret the confusion. The
 manifolds $\mathcal{N}_b$ (indeed all metrics \eqref{eq:walker-metric}) are pp-waves, but not plane-waves since $f_{yy}$ is not quadratic
 in $y$.
\end{remark}

Note that the Cahen-Wallach symmetric spaces $\mathcal{CW}_\varepsilon$ are geodesically complete (see, for example \cite{CLPTV}). The geodesic equations corresponding to the manifolds $\mathcal{N}_b$ can be integrated explicitly. Indeed, it follows from the expressions above for the covariant derivatives that the geodesic equations become
$$
x''(t)=0, \quad
y''(t)=-x'(t)^2 b^{-1}e^{b y(t)}, \quad
\tilde x''(t)= 2 x'(t) y'(t) b^{-1}e^{b y(t)}.
$$
Put $x(t)=\alpha t + \beta$. If $\alpha=0$, the solutions of the system above are clearly globally defined. For non-zero $\alpha$, the second equation above becomes
$y''(t)=-\alpha^2 b^{-1}e^{b y(t)}$ and its solutions are given by
$$
y(t)=b^{-1}\ln\left(\frac{b^2\, C_1^2}{2\alpha^2}\operatorname{sech}^2\left(\frac{1}{2}\mid\! b\, C_1(t+C_2)\!\mid\right)\right),
$$
for some constants $C_1$, $C_2$. Now it follows that the solutions of the geodesic equations are defined for all $t\in\mathbb{R}$, and hence $\mathcal{N}_b$ is geodesically complete.
On the other hand, the plane waves $\mathcal{P}_c$ are not geodesically complete \cite{BL}.  Moreover, one has

\begin{theorem}\label{thm-1.6}
There is a geodesic $\gamma(t)$ in $\mathcal{P}_c$ which defined for $t\in[0,1)$
and there exists a parallel vector field $Y(t)$ along $\gamma(t)$ with
$$
g_f(\dot\gamma,\dot\gamma)=1,\quad g_f(\dot\gamma,Y)=0,\quad
g_f(Y,Y)=1,\quad
\lim_{t\rightarrow1^-}R(\dot\gamma,Y,Y,\dot\gamma)=\infty\,.
$$
Consequently, $\mathcal{P}_c$ is geodesically incomplete and can not be embedded
in a geodesically complete manifold.
\end{theorem}

\noindent{\normalsize\bf Proof.}
We begin by considering a special case.
Let $M=(-\infty,1)\times\mathbb{R}^2$ and let $f=(1-x)^{-2}y^2$.
Let $\gamma(t)=(t,\phi(t),\psi(t))$. Then
$$
\gamma_*(\partial_t)=\partial_x+\phi_t\partial_y+\psi_t\partial_{\tilde x}\,,
$$
and $\gamma$ is a geodesic if and only if
\begin{eqnarray*}
0&=&\{\phi_{tt}+2(1-t)^{-2}\phi\}\partial_y
+\{\psi_{tt}-2(1-t)^{-3}\phi^2-4\phi_t(1-t)^{-2}\phi\}\partial_{\tilde x}\,.
\end{eqnarray*}
Thus we must solve the equations:
\begin{eqnarray}
&&\phi_{tt}+2(1-t)^{-2}\phi=0,\label{eqn-1.d}\\
&&\psi_{tt}-2(1-t)^{-3}\phi^2-4\phi_t(1-t)^{-2}\phi=0\,.\label{eqn-1.e}
\end{eqnarray}
Consider a complex solution $\Phi(t)=(1-t)^\lambda$ for $\lambda\in\mathbb{C}$. To ensure that $\Phi$
solves Equation~(\ref{eqn-1.d}), we require that:
$$
0=\{\lambda(\lambda-1)+2\}(1-t)^{\lambda-2}\quad\text{so}\quad
     \lambda_\pm:=\frac{1\pm\sqrt{1-8}}2\,.
$$
The \emph{real} solutions of (\ref{eqn-1.d}) are then are of the form
\begin{eqnarray*}
\phi(t)&=&\alpha_1e^{\ln(1-t)(1+\sqrt{-1}\sqrt{7})/2}
              +\alpha_2e^{\ln(1-t)(1-\sqrt{-1}\sqrt{7})/2}\\
&=&(1-t)^{\frac12}\{a_1\cos(\sqrt{7}\ln(1-t)/2)+a_2\sin(\sqrt{7}\ln(1-t)/2)\}
\end{eqnarray*}
for suitably chosen constants $(\alpha_1,\alpha_2)$ or, equivalently, for suitably chosen
$(a_1,a_2)$. We set $a_1=1$ and $a_2=1+(\sqrt{7})^{-1}$
to ensure that $\phi(0)=1$ and $\phi_t(0)=0$. We integrate twice to
solve Equation~\eqref{eqn-1.e}
with initial conditions $\psi(0)=0$ and $\psi_t(0)=3$ to ensure:
$$
\gamma(0)=(0,1,0), \quad\gamma_*(\partial_t)(0)=\partial_x+3\partial_{\tilde x},\quad
g_f(\gamma_*\partial_t,\gamma_*\partial_t)(0)=1\,.
$$
Let $Y(t)=\partial_y+\Psi(t)\partial_{\tilde x}$ be a vector field along $\sigma$;
we choose $\Psi$ so that $\Psi(0)=0$ and so that $Y$ is parallel; $\Psi$ satisfies
the partial differential equation
$$\Psi_t-2(1-t)^{-2}\phi=0\,.$$
The plane $\sigma:=\operatorname{Span}\{\dot\gamma,Y\}$ is then
a spacelike $2$-plane and $\{\dot\gamma,Y\}$ is an orthonormal basis for $\sigma$.
Thus the sectional curvature is given by:
$$
R(\dot\gamma,\partial_y,\partial_y,\dot\gamma)=f_{yy}=2(1-t)^{-2}.
$$
Theorem~\ref{thm-1.6} now follows in this special case.

For more general $C$ where $f=C(1-t)^{-2}y^2$,
then we get the corresponding eigenvalue equation
$0=\lambda(\lambda-1)+2C$. The argument will be
essentially the same except that we will sometimes get real solutions for
$\lambda$ depending on the value of
$C$. But the point will be that we can always
solve the equations
 for $t\in[0,1)$ and, by choosing
$\psi$ appropriately, get a spacelike $2$-plane. A similar
argument holds if we take $f=C(t+1)^{-2}y^2$ and let $\gamma(t)=(-t,\star,\star)$.
\hfill$\square$
\medskip

\subsection{Geometric solitons}\label{se:gs}

The objective of the different geometric evolution equations  is to improve a given initial metric by
considering a   flow associated to the geometric object under consideration. The Ricci, Yamabe, and
mean curvature flows are examples extensively studied in the literature. Under suitable conditions,
the Ricci flow evolves an initial metric to an Einstein metric while the Yamabe flow evolves an
initial metric to a new one with constant scalar curvature within the same conformal class. There
are however certain metrics which, instead of evolving by the flow, remain invariant up to scaling and diffeomorphisms, i.e., they are self-similar solutions of the flow.

For any solution of the form $g(t)=\sigma(t)\psi_t^* g(0)$, where $\sigma(t)$ is a smooth function and
$\{\psi_t\}$ a one-parameter family of diffeomorphisms of $M$, there exists a vector field $X$ (the soliton vector field) which relates the Lie derivative of the metric $\mathcal{L}_Xg$ with the geometric object defining the flow under consideration.

\subsubsection{Ricci solitons}\label{se:rs}
A \emph{Ricci soliton} is a pseudo-Riemannian manifold $(M,g)$ which admits
a smooth vector field $X$ (which is called a soliton vector field) on $M$ such that
\begin{equation}\label{soliton}
{\mathcal L}_{X}g+Ric=\lambda g,
\end{equation}
where ${\mathcal L}_{X}$ denotes the Lie derivative in the direction of $X$,
$Ric$ is the Ricci tensor, and $\lambda$ is a real number
($\lambda=\frac{1}{n}(2\text{div} X+Sc)$, where $n=\dim\,M$ and $Sc$ denotes the scalar curvature of $(M,g)$).
A Ricci soliton is said to be \emph{shrinking, steady} or \emph{expanding},
 if $\lambda >0,$ $\lambda =0$ or $\lambda <0$, respectively. Moreover we say
that a Ricci soliton $(M,g)$ is a \emph{gradient Ricci soliton} if the
vector field $X$ satisfies $X=\text{grad}\, h$, for some potential function $h$. In such a case Equation \eqref{soliton} can be written in terms of $h$ as
$2\text{Hes}_{h}+Ric=\lambda g$.
%\begin{equation}\label{gradsoliton}
%2\text{Hes}_{h}+Ric=\lambda g.
%\end{equation}

Three-dimensional Walker metrics admitting a non-trivial (i.e., not Einstein) gradient Ricci soliton were completely described in \cite{BVGRGF}, where it is shown that one of the following two possibilities must occur
\begin{enumerate}
\item[(R.1)] There exist coordinates $(x,y,\tilde x)$ so that the metric $g$ takes the form
\eqref{eq:walker-metric} for some function $f$ satisfying $f_{yy}^{-1}f_{yyy}=b$. Hence
$$
f(x,y)=\frac{1}{\kappa^2}e^{\kappa y}\alpha(x)+y \beta(x)+\gamma(x)
$$
for some arbitrary functions $\alpha(x)$, $\beta(x)$ and $\gamma(x)$.
The potential function of the soliton is given by
$$
h(x,y,\tilde x)=\frac{\kappa}{2} y + \hat{h}(x),\,\, \mbox{where}\,\,\,
\hat{h}_{xx}=\frac{\kappa}{2}\beta(x),
$$
and the soliton vector field is spacelike and given by
$\operatorname{grad} h=\frac{\kappa}{2}\partial_y+\hat{h}_x(x)\partial_{\tilde x}$.

\item[(R.2)] There exist coordinates $(x,y,\tilde x)$ so that the metric $g$ takes the form
\eqref{eq:walker-metric} for some function $f$ satisfying $f_{yyy}=0$. Hence
$$
f(x,y)=y^2 \alpha(x)+y \beta(x)+ \gamma(x)
$$
for some arbitrary functions $\alpha(x)$, $\beta(x)$ and $\gamma(x)$. The potential function of the soliton is given by
$$
h(x,y,\tilde x)=\hat{h}(x),\,\, \mbox{where}\,\,\, \hat{h}_{xx}=- \alpha(y),
$$
and the soliton vector field is lightlike and given by
$\operatorname{grad} h=\hat{h}_x(x)\partial_{\tilde x}$.
\end{enumerate}
Moreover, in both cases the Ricci soliton is steady.
As we shall see in Section \ref{sect-1.3} all gradient Ricci solitons above are $1$-curvature homogeneous, provided that $f_{yy}$ has constant sign.

Further, note that all metrics corresponding to the second case above are plane waves (since the function $f$ is quadratic on $y$), and hence they admit non-gradient vector fields $X$ resulting in expanding and shrinking Ricci solitons \cite{BBGG}. However metrics corresponding to the first case above only admit steady Ricci solitons

\begin{theorem}\label{thm-2.7}
A gradient Ricci soliton $\mathcal{M}_f$ admits a vector field $X$ resulting in a non-steady Ricci soliton if and only if it is a locally conformally flat metric.
\end{theorem}

\noindent{\normalsize\bf Proof.}
Recall here that two Ricci soliton vector fields differ in a homothetic vector field \cite{BBGG}. Hence, a gradient Ricci soliton Walker metric admits a vector field resulting in a non-steady Ricci soliton if and only if it admits a non-Killing homothetic vector field.

Locally conformally flat Walker metrics are plane waves \cite{CGRVA}, and thus they admit homothetic vector fields resulting in expanding, steady and shrinking Ricci solitons \cite{BBGG}.
In the non-locally conformally flat case (i.e., $f_{yyy}\neq 0$), since any gradient Ricci soliton Walker metric satisfies $f_{yyy}=b f_{yy}$, it also follows that $f_{yyyy}=\kappa f_{yy}$ for some $\kappa\neq 0$, and the result follows from Theorem~\ref{th-1}
\hfill$\square$

\subsubsection{Cotton solitons}\label{se:cs}

The {\it Schouten tensor} of any pseudo-Riemannian manifold is given by $\displaystyle
S_{ij}=Ric_{ij}-\frac{Sc}{4}g_{ij}$. Then the \emph{Cotton tensor},
$\displaystyle
C_{ijk}=(\nabla_i S)_{jk}-(\nabla_j S)_{ik}
$
measures the failure of the Schouten tensor  to be a Codazzi tensor. The Cotton tensor is the unique conformal invariant in  dimension three
and it vanishes if and only if the manifold is locally conformally flat. Using the Hodge
$\star$-operator, the $(0,2)$-Cotton tensor is given by
$\displaystyle
C_{ij}=\frac{1}{2\sqrt{g}}C_{nmi}\epsilon^{nm\ell}g_{\ell j},
$
where $\epsilon^{123}=1$.  Moreover,  the
$(0,2)$-Cotton tensor is trace-free and divergence-free \cite{York71} and it  appears naturally in
many physical contexts (see, for example
\cite{Chow-Pope-Sezgin10}, \cite{Garcia-Hehl-Heinicke04} and the references therein).

The only non-zero component of the $(0,2)$-Cotton tensor of a Walker manifold $\mathcal{M}_f$ is given by
$C(\partial_x,\partial_x)=-\frac{1}{2} f_{yyy}$
(and hence the manifold is locally conformally flat if and only if $f_{yyy}=0$).

A geometric flow associated to the Cotton tensor was introduced in \cite{Ali-Ozgur08} as
$$\frac{\partial}{\partial_t}g(t)=-\lambda\, C_{g(t)},
$$
where $C_{g(t)}$ is the  $(0,2)$-Cotton tensor corresponding to $(M,g(t))$. Then one naturally
considers soliton solutions of the Cotton flow. Following  \cite{Ali-Ozgur08}, a \emph{Cotton soliton} is a triple $(M,g,X)$
of a three-dimensional manifold and a vector field $X$ satisfying
\begin{equation}\label{Cotton equation}
\mathcal{L}_X g+C=\lambda\, g,
\end{equation}
where $\lambda$ is a real number.  The Cotton soliton is said to be   \emph{shrinking}, \emph{steady} or \emph{expanding}  if $\lambda>0$, $\lambda = 0$ or $\lambda < 0$, respectively.

The necessary and sufficient conditions for a Walker manifold to be a gradient Cotton soliton were discussed in \cite{CLGRVL}, where it is shown that $\mathcal{M}_f$ is a non-trivial (i.e., not locally conformally flat) gradient Cotton soliton if and only if it is steady and $f_{yyyy}=\kappa f_{yy}$ for some non-zero constant $\kappa$.
Now one of the following three possibilities must occur
\begin{enumerate}
\item[(C.1)] There exist coordinates $(x,y,\tilde x)$ so that the metric $g$ takes the form
\eqref{eq:walker-metric} for some function $f$ satisfying $f_{yyyy}=\kappa^2 f_{yy}$ . Hence
$$
f(x,y)=\frac{1}{\kappa^2}(e^{\kappa y}\alpha_1(x)+e^{-\kappa y}\alpha_2(x))+y \beta(x)+\gamma(x)
$$
where  $\alpha_1(x)$, $\alpha_2(x)$, $\beta(x)$ and $\gamma(x)$ are arbitrary functions.
Moreover, the potential function of the soliton is given by
$h(x,y,\tilde x)=\frac{\kappa}{2} y + \hat{h}(x)$, where
$$
\hat{h}_{xx}(x)=\frac{\kappa}{2}(e^{\kappa y}\alpha_1(x)-e^{-\kappa y}\alpha_2(x)+2\kappa\beta(x)).
$$
The soliton vector field is spacelike and given by
$\operatorname{grad} h=\kappa^2\partial_y+\hat{h}_x(x)\partial_{\tilde x}$.

\item[(C.2)] There exist coordinates $(x,y,\tilde x)$ so that the metric $g$ takes the form
\eqref{eq:walker-metric} for some function $f$ satisfying $f_{yyyy}=-\kappa^2 f_{yy}$ . Hence
$$
f(x,y)=-\frac{1}{\kappa^2}(\cos(\kappa y)\alpha_1(x)+\sin(\kappa y)\alpha_2(x))+y \beta(x)+\gamma(x)
$$
where  $\alpha_1(x)$, $\alpha_2(x)$, $\beta(x)$ and $\gamma(x)$ are arbitrary functions.
Moreover, the potential function of the soliton is given by
$h(x,y,\tilde x)=-\frac{\kappa}{2} y + \hat{h}(x)$, where
$$
\hat{h}_{xx}(x)=\frac{\kappa}{2}(\cos(\kappa y)-\sin(\kappa y) -2\kappa\beta(y)).
$$
The soliton vector field is spacelike and given by
$\operatorname{grad} h=-\kappa^2\partial_y+\hat{h}_x(x)\partial_{\tilde x}$.

\item[(C.3)] There exist coordinates $(x,y,\tilde x)$ so that the metric $g$ takes the form
\eqref{eq:walker-metric} for some function $f$ satisfying $f_{yyyy}=0$ . Hence
$$
f(x,y)=y^3 \alpha_1(x)+y^2\alpha_2(x)+y \beta(x)+\gamma(x)
$$
where $\alpha_1(x)$, $\alpha_2(x)$, $\beta(x)$ and $\gamma(x)$ are arbitrary functions.
Moreover, the potential function of the soliton is given by
$h(x,y,\tilde x)= \hat{h}(x)$, where
$$
\hat{h}_{xx}(x)=-3\alpha_1(x).
$$
The soliton vector field is spacelike and given by
$\operatorname{grad} h=-\kappa^2\partial_y+\hat{h}_x(x)\partial_{\tilde x}$.

\end{enumerate}
Moreover, in all cases above the gradient Cotton soliton is steady.

Note that two Cotton soliton vector fields
$(\mathcal{L}_{X_i}g+C=\lambda_i g, \,\, i=1,2)$ differ by a homothetic vector field since
$$
\mathcal{L}_{X_1-X_2}g=(\lambda_1-\lambda_2)g.
$$
Hence, as well as for Ricci solitons, no Walker metric corresponding to (C.1) and (C.2) supports any non-trivial Cotton soliton of non-steady type as a consequence of the following.

\begin{theorem}\label{th-1}
Let $\mathcal{M}_f$ be a non-flat Walker manifold satisfying $f_{yyyy}=bf_{yy}$ $(b\neq 0)$. Then any homothetic vector field on $\mathcal{M}_f$ is necessarily a Killing vector field.
\end{theorem}

\noindent{\normalsize\bf Proof.}
A vector field $X=\mathcal{A}(x,y,\tilde x)\partial_x+\mathcal{B}(x,y,\tilde x)\partial_y+\mathcal{C}(x,y,\tilde x)\partial_{\tilde x}$ is a homothetic vector field if and only if
$\mathcal{L}_Xg=\mu g$ for some constant $\mu$.
A straightforward calculation shows that $X$ is a homothetic vector field if and only if
\begin{equation}\label{eq:homoth-1}
\begin{array}{lcl}
\mathcal{C}_y-2f\mathcal{A}_y+\mathcal{B}_x=0,& &\mathcal{B}_{\tilde x}+\mathcal{A}_y=0,\\
\noalign{\medskip}
\mathcal{C}_{\tilde x}- 2f\mathcal{A}_{\tilde x}+\mathcal{A}_x=\mu, & &
2\mathcal{B}_y=\mu,\\
\noalign{\medskip}
\mathcal{B}f_y+\mathcal{A}f_x-\mathcal{C}_x-f(\mu-2\mathcal{A}_x)=0, & &\mathcal{A}_{\tilde x}=0,\\
\end{array}
\end{equation}
Now a standard integration process shows that all solutions of \eqref{eq:homoth-1} must take the form
$$
X=\{ a x+\overline{a}\}\partial_x+\left(\frac{\mu}{2}y+U(x)\right)\partial_y+
\{T(x)-yU_x(x)+\tilde x(\mu-\overline{a})\}\partial_{\tilde x},
$$
for some constants $a$, $\overline{a}$ and some functions $U(x)$, $T(x)$,
where
\begin{equation}\label{eq:homoth-2}
2 a f+\mu +T_x(x)-y U_{xx}(x)-\left(\frac{\mu}{2}y +U(x)\right)f_y-(a x +\overline{a})f_x=0.
\end{equation}
Next, differentiate twice in \eqref{eq:homoth-2} with respect to $y$ to obtain
\begin{equation}\label{eq:homoth-3}
2 a f_{yy}+\left(\frac{\mu}{2}y+U(x)\right)f_{yyy} +(a x +\overline{a})f_{xyy}=0.
\end{equation}
Differentiating once again with respect to $y$ and using that $f_{yyyy}=bf_{yy}$, one gets
\begin{equation}\label{eq:homoth-4}
\left(2 a+\frac{\mu}{2}\right) f_{yyy}+b\left(\frac{\mu}{2}y+U(x)\right)f_{yy} +(a x +\overline{a})f_{xyyy}=0.
\end{equation}
A further differentiation in \eqref{eq:homoth-4} with respect to $y$ (using that $f_{yyyy}=bf_{yy}$) gives
\begin{equation}\label{eq:homoth-5}
b\left(2 a+\frac{\mu}{2}\right) f_{yy}+b\frac{\mu}{2}f_{yy}+b\left(\frac{\mu}{2}y+U(x)\right)f_{yyy} +b(a x +\overline{a})f_{xyy}=0
\end{equation}
which, together with \eqref{eq:homoth-3}, shows that $b\mu f_{yy}=0$. Hence $X$ is Killing.
\hfill$\square$

\begin{remark}\rm
Let $\mathcal{M}_f$ be given by
$
f(x,y)=y^3 e^{-\lambda x}+y^2+y e^{\lambda x}+\gamma(x)$.
Then the vector field
$X=\frac{1}{2}\partial_x+\frac{\lambda}{2}y\partial_y+\left\{\lambda \tilde{x}+\theta(x)\right\}\partial_{\tilde x}$
is a Cotton soliton for any function $\theta(x)$ which satisfies the identity
$\theta_x=3 e^{-\lambda x}+(\frac{1}{2}-\lambda)\gamma(x)$.
Moreover the Cotton soliton is expanding or shrinking depending on the sign of $\lambda$.
Hence there are gradient Cotton solitons (C.3) which also admit non-Killing homothetic vector fields.
\end{remark}

\subsection{Curvature homogeneity}

A pseudo-Riemannian manifold $(M,g)$ is said to be \emph{$k$-curvature homogeneous} if
for each pair of points $p,q\in M$ there is a linear isometry
$\Phi_{pq}:T_pM\rightarrow T_qM$ such that
\[
\Phi_{pq}^*R(q)=R(p), \,\, \Phi_{pq}^*\nabla R(q)=\nabla R(p), \,\,
\dots, \,\, \Phi_{pq}^*\nabla^k R(q)=\nabla^k R(p)
\]
where $R$, $\nabla R$, $\dots$, $\nabla^k R$ stands for the
curvature tensor and its covariant derivatives up to order $k$.

Clearly any locally homogeneous manifold is curvature homogeneous
and the converse holds true if $k$ is sufficiently large.
An open question in the study of curvature homogeneity is to
decide the minimum level of curvature homogeneity needed to show that
a space is locally homogeneous. A general estimate of the form $k_M+1\leq n(n-1)/2$
(where $n$ is the dimension of the manifold) was obtained by Singer \cite{Singer} (see also \cite{PS}).
However, there are sharper bounds in low dimensions.
A Riemannian manifold which is $1$-curvature homogeneous is locally homogeneous in dimension $\leq 4$.
However one needs $2$-curvature homogenenity to ensure local homogeneity in the three-dimensional Lorentzian
setting (see \cite{G07} for more information and references).

A $0$-curvature homogeneous manifold is said to be modeled on a symmetric space if its curvature tensor at each point is that of a symmetric space.
A complete and simply connected indecomposable Lorentzian symmetric space is either irreducible, and hence of constant sectional curvature, or otherwise it is a Cahen-Wallach symmetric space \cite{CLPTV}.
Now, an immediate application of Schur's lemma shows that {a curvature homogeneous Lorentz manifold
modeled on an irreducible symmetric space has constant sectional curvature}.
On the other hand, curvature homogeneous Lorentzian manifolds modeled on indecomposable symmetric spaces need not to be symmetric, but they are Walker manifolds \cite{CLPTV}.

\subsection{Summary of results}

The purpose of this work is to analyze the class of three-dimensional manifolds with recurrent curvature under different curvature homogeneity assumptions. The question of $1$-curvature homogeneity is dealt with by the following result  which will be proved in Section~\ref{sect-1.4}:

 \begin{theorem}\label{thm-1.3}
 Assume that $f_{yy}>0$ and that $f_{yy}$ is non-constant.
Then $\mathcal{M}_f$ is $1$-curvature homogeneous if
and only if exactly one of the following two possibilities holds:
\begin{enumerate}
\item $f_{yy}(x,y)=\alpha(x)e^{by}$ where
$0\ne b\in\mathbb{R}$ and where
$\alpha(x)$ is arbitrary. This manifold is $1$-curvature modeled on the manifold $\mathcal{N}_b$
of Definition~\ref{def-2.3}.
\item $f_{yy}(x,y)=\alpha(x)$ where $\alpha=c\cdot\alpha_x^{3/2}$ for
some $0\ne c\in\mathbb{R}$. This manifold is locally homogeneous and
is locally isometric to the manifold $\mathcal{P}_c$ of Definition~\ref{def-2.3}.
\end{enumerate}
 \end{theorem}

Previous result coupled with those in Section~\ref{se:rs} show that

\begin{theorem}\label{thm-2.9}
Any $1$-curvature homogeneous Lorentzian three-manifold with recurrent curvature is a steady gradient Ricci soliton.
\end{theorem}

The complete description of all $2$-curvature homogenenous $\mathcal{M}_f$ will be discussed
 in Sections ~\ref{sect-1.5x} and \ref{sect-1.6} as well as in Section~\ref{sect-1.7}.
 The following result shows that any locally homogeneous Walker three-dimensional manifold is locally isometric to one of the models introduced in Definition \ref{def-2.3}. (This clarifies the results in \cite{BCDL}).

 \begin{theorem}\label{thm-1.4}
 The manifold $\mathcal{M}_f$ is $2$-curvature homogeneous if and only if
 it falls into one of the three families:
 \begin{enumerate}
 \item $f=b^{-2}\alpha(x)e^{by}+\beta(x)y+\gamma(x)$  for
 $\beta(x)=b^{-1}\alpha^{-1}\{\alpha_{xx}-\alpha_x^2\alpha^{-1}\}$
 where $b\ne0$ and $\alpha>0$.
 The manifold $\mathcal{M}_f$ is locally isometric to the manifold
 $\mathcal{N}_{b}$ of Definition~\ref{def-2.3}
 and consequently is locally homogeneous.
 \smallbreak\item $f_{yy}=\alpha(x)>0$ where $\alpha_x=c\alpha^{3/2}$
  for $c>0$; thus $\alpha=\tilde c(x-x_0)^{-2}$ for some $(\tilde c,x_0)$.
 The manifold $\mathcal{M}_f$ is locally isometric to the manifold
 $\mathcal{P}_c$ of Definition~\ref{def-2.3} and consequently locally homogeneous.
\smallbreak \item $f=\varepsilon y^2+\beta(x)y+\gamma(x)$
 where $0<\varepsilon\in\mathbb{R}$. The manifold $\mathcal{M}_f$
 is locally isometric
 to the manifold $\mathcal{CW}_\varepsilon$ of Definition~\ref{def-2.3} and consequently is locally homogeneous.
 \end{enumerate}
 \end{theorem}

Hence, as an application of the results in Section \ref{se:cs}
\begin{theorem}\label{thm-2.10}
Any $2$-curvature homogeneous Lorentzian three-manifold with recurrent curvature is a steady gradient Cotton soliton.
\end{theorem}

Finally note that the model manifolds $\mathcal{P}_c$ and $\mathcal{CW}_\varepsilon$ admit Ricci and Cotton solitons of any kind (expanding, steady and shrinking), but $\mathcal{N}_b$ only admits steady Ricci and Cotton solitons.

\section{Higher order curvature homogeneity}\label{section-3}

We introduce the following $1$-curvature and $2$-curvature models which will play an important
role in our development:

\begin{definition}\label{defn-x}\rm
Let $V=\operatorname{span} \{\xi_1,\xi_2,\xi_3\}$ be a $3$-dimensional vector space.
Let $\mathfrak{A}_0$ be the $0$-curvature model defined by:
$$
\langle\xi_1,\xi_3\rangle=\langle\xi_2,\xi_2\rangle=1,
\quad
\mathcal{R}(\xi_1,\xi_2,\xi_2,\xi_1)=1.
$$
\begin{enumerate}
%\item Let $\mathfrak{N}_0$ be the $0$-curvature model defined by:
%$$
%\langle\xi_1,\xi_3\rangle=\langle\xi_2,\xi_2\rangle=1,
%\quad
%\mathcal{R}(\xi_1,\xi_2,\xi_2,\xi_1)=1.
%$$
\item
\begin{enumerate}
\item Let $\mathfrak{N}_1(b)$ be the $1$-curvature model induced from
$\mathfrak{A}_0$ by
imposing a condition on $\nabla\mathcal{R}$:
$$
\begin{array}{c}
\langle\xi_1,\xi_3\rangle=\langle\xi_2,\xi_2\rangle=1,
\,\,
\mathcal{R}(\xi_1,\xi_2,\xi_2,\xi_1)=1,\\
\noalign{\medskip}
\nabla\mathcal{R}(\xi_1,\xi_2,\xi_2,\xi_1;\xi_1)=0,
\qquad\qquad
\nabla\mathcal{R}(\xi_1,\xi_2,\xi_2,\xi_2;\xi_2)=b,
\end{array}
$$
for some non-zero $b\in\mathbb{R}$.
\item Let $\mathfrak{P}_1(c)$ be the $1$-curvature model induced from
$\mathfrak{A}_0$ by
$$
\begin{array}{c}
\langle\xi_1,\xi_3\rangle=\langle\xi_2\xi_2\rangle=1,
\,\,
\mathcal{R}(\xi_1,\xi_2,\xi_2,\xi_1)=1,\\
\noalign{\medskip}
\nabla\mathcal{R}(\xi_1,\xi_2,\xi_2,\xi_1;\xi_1)=c,
\qquad\qquad
\nabla\mathcal{R}(\xi_1,\xi_2,\xi_2,\xi_2;\xi_2)=0,
\end{array}
$$
for some non-zero $c\in\mathbb{R}$.
\end{enumerate}

\item
\begin{enumerate}
\item Let $\mathfrak{N}_{2}(b)$ be the $2$-curvature model induced from
$\mathfrak{N}_1(b)$ by
imposing a condition on $\nabla^2\mathcal{R}$:
$$
\begin{array}{c}
\langle\xi_1,\xi_3\rangle=\langle\xi_2,\xi_2\rangle=1,
\,\,
\mathcal{R}(\xi_1,\xi_2,\xi_2,\xi_1)=1,\\
\noalign{\medskip}
\begin{array}{ll}
\nabla\mathcal{R}(\xi_1,\xi_2,\xi_2,\xi_1;\xi_1)=0,
&
\nabla\mathcal{R}(\xi_1,\xi_2,\xi_2,\xi_1;\xi_2)=b,\\
\noalign{\medskip}
\nabla^2\mathcal{R}(\xi_1,\xi_2,\xi_2,\xi_1;\xi_1,\xi_1)=-1,
&
\nabla^2\mathcal{R}(\xi_1,\xi_2,\xi_2,\xi_1;\xi_2,\xi_2)=b^2\\
\noalign{\medskip}
\nabla^2\mathcal{R}(\xi_1,\xi_2,\xi_2,\xi_1;\xi_1,\xi_2)=0,
&
\nabla^2\mathcal{R}(\xi_1,\xi_2,\xi_2,\xi_1;\xi_2,\xi_1)=0.
\end{array}
\end{array}
$$
\item Let $\mathfrak{P}_{2}(c)$ be the $2$-curvature model induced from
$\mathfrak{P}_1(c)$ by
imposing a condition on $\nabla^2\mathcal{R}$:
$$
\begin{array}{c}
\langle\xi_1,\xi_3\rangle=\langle\xi_2\xi_2\rangle=1,
\,\,
\mathcal{R}(\xi_1,\xi_2,\xi_2,\xi_1)=1,\\
\noalign{\medskip}
\begin{array}{ll}
\nabla\mathcal{R}(\xi_1,\xi_2,\xi_2,\xi_1;\xi_1)=c,
&
\nabla\mathcal{R}(\xi_1,\xi_2,\xi_2,\xi_1;\xi_2)=0,\\
\noalign{\medskip}
\nabla^2\mathcal{R}(\xi_1,\xi_2,\xi_2,\xi_1;\xi_1,\xi_1)=\frac{3}{2}c^2,
&
\nabla^2\mathcal{R}(\xi_1,\xi_2,\xi_2,\xi_1;\xi_2,\xi_2)=0,\\
\noalign{\medskip}
\nabla^2\mathcal{R}(\xi_1,\xi_2,\xi_2,\xi_1;\xi_1,\xi_2)=0,
&
\nabla^2\mathcal{R}(\xi_1,\xi_2,\xi_2,\xi_1;\xi_2,\xi_1)=0.
\end{array}
\end{array}
$$
\end{enumerate}
\end{enumerate}
\end{definition}

\subsection{$0$-curvature homogeneity}\label{sect-1.3}
We suppose $f_{yy}>0$ henceforth in our study of $0$-curvature homogeneity; the case $f_{yy}<0$ is completely analogous.  We have
$$
\ker(\mathcal{R})=\operatorname{Span}\{\partial_{\tilde x}\}\quad\text{and}\quad
\operatorname{Range}(\mathcal{R})
=\operatorname{Span}\{\partial_{\tilde x},\partial_y\}\,.
$$
Thus these two subspaces are invariantly defined. We set
$$\xi_1:=a_{11}(\partial_x+f\partial_{\tilde x}+
     a_{12}\partial_y+a_{13}\partial_{\tilde x}),\quad
     \xi_2:=a_{22}\partial_y+a_{23}\partial_{\tilde x},\quad
     \xi_3:=a_{33}\partial_{\tilde x}\,,
$$
to be a pseudo-orthonormal basis. We wish to choose the basis to see that $\mathcal{M}_f$ is 0-curvature modeled on the 0-model $\mathfrak{A}_0$ of Definition \ref{defn-x}. To ensure that the inner product is properly normalized, we need the equations:
$$
2a_{13}+a_{12}^2=0,\quad a_{23}+a_{12}a_{22}=0,\quad
    a_{22}^2=1,\quad a_{11}a_{33}=1\,.
$$
To ensure that $\mathcal{R}(\xi_1,\xi_2,\xi_2,\xi_1)=1$, we require:
$$
a_{11}^2a_{22}^2f_{yy}=1\,.
$$
Let $a_{12}$ be arbitrary for the moment; we will normalize this
parameter subsequently in Section~\ref{sect-1.4}. We have:
\begin{equation}\label{eqn-1.a}
\begin{array}{rrr}
a_{11}=f_{yy}^{-1/2},&a_{12}=\star,&a_{13}=-\textstyle\frac12a_{12}^2,\\
\noalign{\medskip}
a_{22}=1,&a_{23}=-a_{12},&a_{33}=f_{yy}^{1/2}.
\end{array}
\end{equation}
The parameters $a_{13}$, $a_{23}$, and $a_{33}$ play no further role.
This shows that $\mathcal{M}_f$ is in fact 0-curvature modeled on $\mathfrak{A}_0$. Taking $\mathcal{M}_f$ to define $\mathcal{CW}_\varepsilon$, then yields:

\begin{lemma}\label{lem-1.6}
If $f_{yy}>0$, then $\mathcal{M}_f$ is $0$-curvature homogeneous modeled on the symmetric
space $\mathcal{CW}_\varepsilon$.
\end{lemma}

We shall assume henceforth that $f_{yy}$ is non-constant.

\subsection{1-curvature homogeneity -- the proof of Theorem~\ref{thm-1.3}}\label{sect-1.4}

Adopt the normalizations of Equation~(\ref{eqn-1.a})
to ensure $R(\xi_1,\xi_2,\xi_2,\xi_1)=1$. We then have:
$$
\nabla R(\xi_1,\xi_2,\xi_2,\xi_1;\xi_2)=f_{yyy}\cdot f_{yy}^{-1}\,.
$$
Since $a_{12}$ plays no role in $\nabla R(\xi_1,\xi_2,\xi_2,\xi_1;\xi_2)$, we see
that $f_{yyy}\cdot f_{yy}^{-1}$ is an isometry invariant.
Consequently, if $\mathcal{M}_f$ is
1-curvature homogeneous, then $f_{yyy}=b\cdot f_{yy}$
for some $b\in\mathbb{R}$ and
thus $f_{yyy}=\alpha(x)e^{by}$. The possibility in Assertion (1) arises from $b\ne0$
and the possibility in Assertion (2) arises from $b=0$ and $\alpha(x)$ non-constant.

If $b\ne0$, then $f_{yyy}\ne0$ and we may set:
\begin{equation}\label{eqn-1.b}
\begin{array}{rrr}
a_{11}=f_{yy}^{-1/2},&a_{12}=-f_{xyy}\cdot f_{yyy}^{-1},&a_{13}=-\textstyle\frac12a_{12}^2,\\
\noalign{\medskip}
a_{22}=1,&a_{23}=-a_{12},&a_{33}=f_{yy}^{1/2}.
\vphantom{\vrule height 12pt}
\end{array}\end{equation}
With these normalizations, we establish Assertion (1) by computing:
\begin{eqnarray}
&&R(\xi_1,\xi_2,\xi_2,\xi_1)=f_{yy}\cdot f_{yy}^{-1}=1,\nonumber\\
\noalign{\medskip}
&&\nabla R(\xi_1,\xi_2,\xi_2,\xi_1;\xi_1)=\{f_{xyy}+a_{12}f_{yyy}\}f_{yy}^{-3/2}=0,
    \label{eqn-1.c}\\
\noalign{\medskip}
&&\nabla R(\xi_1,\xi_2,\xi_2,\xi_1;\xi_2)=f_{yyy}\cdot f_{yy}^{-1}=b\,,\nonumber
\end{eqnarray}
thus showing that $\mathcal{M}_f$ is $1$-curvature homogeneous modeled on  the 1-curvature model $\mathfrak{N}_1(b)$ of Definition \ref{defn-x}.

On the other hand, if $b=0$, then $a_{12}$ plays no role in the computation of $\nabla R$ and we have:
\begin{eqnarray*}
&&R(\xi_1,\xi_2,\xi_2,\xi_1)=f_{yy}\cdot f_{yy}^{-1}=1,\\
\noalign{\medskip}
&&\nabla R(\xi_1,\xi_2,\xi_2,\xi_1;\xi_1)=f_{xyy}\cdot f_{yy}^{-3/2}
=\alpha_x\cdot\alpha^{-3/2},\\
\noalign{\medskip}
&&\nabla R(\xi_1,\xi_2,\xi_2,\xi_1;\xi_2)=f_{yyy}\cdot f_{yy}^{-1}=0\,.
\end{eqnarray*}
So $\mathcal{M}_f$ will be $1$-curvature homogeneous if and only if
$\alpha_x=c\cdot\alpha^{3/2}$ for some $0\ne c\in\mathbb{R}$,  thus showing that
$\mathcal{M}_f$ is $1$-curvature homogeneous modeled on  the 1-curvature model $\mathfrak{P}_1(c)$ of Definition \ref{defn-x}.
This establishes the first part of Assertion (2). We postpone the proof of the second part of Assertion (2) to Section \ref{sect-1.7}.

\begin{remark}\rm
We have assumed that $f_{yy}>0$ and thus $\alpha(x)>0$. If we take
$\alpha(x)=a(x-x_0)^{-2}$ for $a>0$ and $x>x_0$, we then have
$$
\alpha_x=-2a(x-x_0)^{-3}=-2a^{-1/2}\{a(x-x_0)^{-2}\}^{3/2}\,.
$$
On the other hand, if $x<x_0$, then
$$
\alpha_x=-2a(x-x_0)^{-3}=2a^{-1/2}\{a(x-x_0)^{-2}\}^{3/2}\,.
$$
So we can get both positive and negative proportionality constants.
And the value at $x=x_1$ can be adjusted by choosing $x_0$ appropriately.
Thus this is the most general possible solution.
\end{remark}

\subsection{2-curvature homogeneity -- Case 1}\label{sect-1.5x}
We continue the discussion of Section~\ref{sect-1.4} and assume
 $f$ has the
form given in Assertion~(1) of Theorem~\ref{thm-1.4}, i.e.
$$
f=b^{-2}\alpha(x)e^{by}+\beta(x)y+\gamma(x)\,.
$$

\begin{lemma}
Given $0\ne b\in\mathbb{R}$ and $\alpha(x)>0$,
set $\beta(x)=b^{-1}\alpha^{-1}\{\alpha_{xx}-\alpha_x^2\alpha^{-1}\}$
so
$$
f=b^{-2}\alpha(x)e^{by}+b^{-1}\alpha^{-1}\{\alpha_{xx}-\alpha_x^2\alpha^{-1}\}y+\gamma(x)\,.
$$
Then $\mathcal{M}_f$ is a $2$-curvature homogeneous manifold which is
$2$-curvature modeled on $\mathcal{N}_b$ for any $\gamma(x)$.
\end{lemma}

\noindent{\normalsize\bf Proof.}
Let $f=\gamma(x)+\beta(x)y+b^{-2}\alpha(x)e^{by}$.
We adopt the normalizations of Equation~(\ref{eqn-1.b}) and continue
the computations of Equation~(\ref{eqn-1.c}) to see:
\begin{eqnarray*}
\nabla^2R(\xi_1,\xi_2,\xi_2,\xi_1;\xi_1,\xi_2)&=&f_{yy}^{-3/2}
\left\{f_{xyyy}-f_{yyyy}f_{xyy}/f_{yyy}\right\}\\
\noalign{\medskip}
&=&f_{yy}^{-3/2}\left\{b\alpha_x-b\alpha_x\right\}
     e^{by}=0,\\
     \noalign{\medskip}
\nabla^2R(\xi_1,\xi_2,\xi_2,\xi_1;\xi_2,\xi_2)&=&f_{yy}^{-1}f_{yyyy}=b^2\,.
\end{eqnarray*}
Thus only $\nabla^2R(\xi_1,\xi_2,\xi_2,\xi_1;\xi_1,\xi_1)$ is
relevant to our discussion.
The term which involves the expression $f_yf_{yyy}=(\beta +b^{-1}\alpha e^{by})b\alpha e^{by}$ is crucial.
We expand:
\begin{eqnarray*}
&&\nabla^2R(\xi_1,\xi_2,\xi_2,\xi_1;\xi_1,\xi_1)=
    f_{yy}^{-2}\{
      \nabla^2R(\partial_x,\partial_y,\partial_y,\partial_x;\partial_x,\partial_x)\\
      \noalign{\smallskip}
&&\qquad+2a_{12}
\nabla^2R(\partial_x,\partial_y,\partial_y,\partial_x;\partial_x,\partial_y)
  +a_{12}^2
  \nabla^2R(\partial_x,\partial_y,\partial_y,\partial_x;\partial_y,\partial_y)\}\\
  \noalign{\medskip}
&&\quad=f_{yy}^{-2}\{f_{xxyy}-f_yf_{yyy}-2f_{xyy}f_{yyy}^{-1}f_{xyyy}
       +f_{xyy}^2f_{yyy}^{-2}f_{yyyy}\}\\
       \noalign{\medskip}
&&\quad=e^{-by}\alpha^{-2}
 \{\alpha_{xx}-b\alpha\beta(x)-2\alpha_x^2\alpha^{-1}+\alpha_x^2\alpha^{-1}\}
 -1\,.
 \end{eqnarray*}
We complete the proof by setting
$$
\beta(x)=b^{-1}\alpha^{-1}\{\alpha_{xx}-\alpha_x^2\alpha^{-1}\}.
$$
The above shows that $\mathcal{M}_f$ given by Theorem \ref{thm-1.4}-(1) is 2-curvature homogeneous (and hence locally homogeneous) modeled on the curvature model $\mathfrak{N}_2(b)$ of Definition \ref{defn-x}.
\hfill$\square$

\subsection{2-curvature homogeneity -- Case 2}\label{sect-1.6}

We complete the proof of Theorem~\ref{thm-1.2} by showing:

\begin{lemma}\label{lem-1.10}
Suppose that $f_{yy}=\alpha(x)>0$ where $\alpha_x=c\alpha^{3/2}$
for $0\ne c\in\mathbb{R}$. Then $\mathcal{M}_f$ is locally homogeneous,
\end{lemma}

\noindent{\normalsize\bf Proof.}
We adopt the normalizations of Equation~(\ref{eqn-1.a}); the parameter
$a_{12}$ plays no role. We compute:
\begin{eqnarray*}
\nabla^2R(\xi_1,\xi_2,\xi_2,\xi_1;\xi_1,\xi_2)&=&\alpha ^{-3/2}f_{xyyy}=0,\\
\noalign{\medskip}
\nabla^2R(\xi_1,\xi_2,\xi_2,\xi_1;\xi_2,\xi_2)&=&\alpha^{-1}f_{yyyy}=0\,.
\end{eqnarray*}
Thus only $\nabla^2R(\xi_1,\xi_2,\xi_2,\xi_1;\xi_1,\xi_1)$ is relevant.
The $f_{yyy}$ term
no  longer plays a role so we have
$$\nabla^2R(\xi_1,\xi_2,\xi_2,\xi_1;\xi_1,\xi_1)=\alpha^{-2}\alpha_{xx}\,.$$
We have $\alpha_x=c\alpha^{3/2}$ and thus
$\alpha_{xx}=\frac{3}{2}c\cdot\alpha_x\cdot\alpha^{1/2}=\frac{3}{2}c^2\alpha^2$, from where it follows that
$\nabla^2R(\xi_1,\xi_2,\xi_2,\xi_1;\xi_1,\xi_1)$ is constant.
Hence $\mathcal{M}_f$ given by Theorem  \ref{thm-1.4}-(2) is 2-curvature homogeneous  modeled on the curvature model $\mathfrak{P}_2(c)$ of Definition \ref{defn-x}.

We proceed
inductively to show that the only non-zero entry in the $k$-th covariant derivative $\nabla^kR$ is
given by $\nabla^kR(\xi_1,\xi_2,\xi_2,\xi_1;\xi_1,...,\xi_1)$ and
that
$$\nabla^{k-1}R(\partial_x,\partial_y,\partial_y,\partial_x;\partial_x,...,\partial_x)
=c_{k-1}\alpha^{(1+k)/2}\,.$$
It then follows that
$$
\begin{array}{rcl}
\nabla^kR(\partial_x,\partial_y,\partial_y,\partial_x;\partial_x,...,\partial_x)
&=&c_{k-1}{\textstyle\frac{1+k}2}\alpha_x\alpha^{(-1+k)/2}\\
\noalign{\medskip}
&=&
c_{k-1}c{\textstyle\frac{1+k}2}\alpha^{(3-1+k)/2}\\
\noalign{\medskip}
&=&c_k\alpha^{(2+k)/2}
\end{array}
$$
for $c_k:=c_{k-1}c{\textstyle\frac{1+k}2}$. It then follows that
$$\nabla^kR(\xi_1,\xi_2,\xi_2,\xi_1;\xi_1,...,\xi_1)=
\alpha^{(-2-k)/2}\nabla^kR(\partial_x,\partial_y,\partial_y;\partial_x,...,\partial_x)
=c_k\,.$$
Thus $\mathcal{M}_f$ is $k$-curvature homogeneous for all $k$ and hence
locally homogeneous.
\hfill$\square$

\section{The isometry classes}\label{sect-1.7}

%\subsection{Change of coordinates}\label{sect-1.7}
We now model the transformations made in Section~\ref{sect-1.3} on the geometric level.
Let $\phi$ and $\psi$ be smooth functions of $x$. We consider the coordinate transformation:
$$
T(x,y,\tilde x)=(x,y+\phi,\tilde x-\phi_xy+\psi)\,.
$$
We then have
$$
T_*\partial_x=\partial_x+\phi_x\partial_y+(-\phi_{xx}y+\psi_x)\partial_{\tilde x},
\quad
T_*\partial_y=\partial_y-\phi_x\partial_{\tilde x},
\quad
T_*\partial_{\tilde x}=\partial_{\tilde x},
$$
and thus
\begin{eqnarray*}
&&g_f(T_*\partial_x,T_*\partial_x)=-2\left\{f(x,y+\phi)
       +\phi_{xx}y-\psi_x-{\textstyle\frac12}\phi_x^2\right\},\\
       \noalign{\medskip}
&&g_f(T_*\partial_x,T_*\partial_{\tilde x})=1,\\
\noalign{\medskip}
&&g_f(T_*\partial_y,T_*\partial_y)=1,\\
\noalign{\medskip}
&&g_f(T_*\partial_x,T_*\partial_y)=-\phi_x+\phi_x=0,\\
\noalign{\medskip}
&&g_f(T_*\partial_y,T_*\partial_{\tilde x})=g_f(T_*\partial_{\tilde x},T_*\partial_{\tilde x})=0\,.
\end{eqnarray*}
Consequently $T^*(g_f)=g_{\tilde f}$, hence defining an isometry between the metrics $g_f$ and
$g_{\tilde f}$, where
$$
\tilde f=f(x,y+\phi)+\phi_{xx}y-\psi_x-{\textstyle\frac12}\phi_x^2\,.
$$
Now we complete the proof of Theorem~\ref{thm-1.4} showing that any locally homogeneous manifold
$\mathcal{M}_f$ is locally isometric to one of $\mathcal{N}_b$, $\mathcal{P}_c$ or $\mathcal{CW}_\varepsilon$ given at Definition \ref{def-2.3}.

\subsubsection{The model manifold $\mathcal{N}_b$ of Definition~\ref{def-2.3}}\label{sect-4.0.1}
Suppose that $f=b^{-2}e^{by}$. Let $\alpha=\alpha(x)>0$ be given and let
$\gamma=\gamma(x)$ be given.
Set $\phi=\ln(\alpha)b^{-1}$. Choose $\psi$ so that
$-\psi_x-\frac12\phi_x^2=\gamma$.
We compute:
\begin{eqnarray*}
\tilde f(x,y)&=&b^{-2}e^{b(y+b^{-1}\ln(\alpha))}
+b^{-1}\{\alpha^{-1}\alpha_x\}_xy+\gamma\\
\noalign{\medskip}
&=&b^{-2}\alpha e^{by}+b^{-1}\alpha^{-1}\{\alpha_{xx}
-\alpha_x^2\alpha^{-1}\}y+\gamma\,.
\end{eqnarray*}
This completes the proof of Theorem~\ref{thm-1.4}-(1).

\begin{remark}\rm
We can use this formalism to show that $\mathcal{N}_b$ is a homogeneous
space. Suppose given a point $(a_1,a_2,a_3)\in\mathbb{R}^3$.
We consider the map:
$$T(x,y,\tilde x)=(e^{-ba_2/2}x+a_1,y+a_2,e^{ba_2/2}\tilde x+a_3)\,.$$
The only inner product which has been changed is
$g_f(T_*\partial_x,T_*\partial_x)=-2e^{by}e^{ba_2}e^{-ba_2}$ and thus
$T^*g_f=g_f$. Thus $g_f$ is an isometry and $T(0,0,0)=(a_1,a_2,a_3)$. Consequently, $\mathcal{N}_{b}$ is a homogeneous space.
\end{remark}

\subsubsection{The manifold $\mathcal{P}_c$ of Definition~\ref{def-2.3}}
Assume that $f=\frac12y^2\alpha_c$. Suppose that $\beta=\beta(x)$ and
$\gamma=\gamma(x)$ are given. Choose $\phi$ so
$\alpha_c\phi+\phi_{xx}=\beta$ and choose $\psi$
so that $-\psi_x-\frac12\phi_x^2+\frac12\alpha_c\phi=\alpha$.
Then:
\begin{eqnarray*}
\tilde f(x,y)&=&{\textstyle\frac12}y^2\alpha_c(x)+y(\alpha_c\phi+\phi_{xx})
-\psi_x-{\textstyle\frac12}\phi_x^2+{\textstyle\frac12}\alpha_c\phi^2\\
\noalign{\medskip}
&=&{\textstyle\frac12}y^2\alpha_c(x)+y\beta+\gamma\,.
\end{eqnarray*}
This completes the proof of Theorem~\ref{thm-1.4}-(2).

\begin{remark}\rm
We can use this formalism to show that $\mathcal{N}_c$ is a local homogeneous
space. We take $f(x,y)=cy^2(x+1)^{-2}$ and $M=(-1,\infty)\times\mathbb{R}^2)$;
the case when $f(x,y)=cy^2(x-1)^{-2}$ and $M=(-\infty,1)$ is handled similarly
(the question of where the singularity is relative to $x=0$ plays an
important role). Let $(a_1,a_2,a_3)\in\mathbb{R}^3$ be given with $a_1>-1$.
 Choose $\phi$ and $\psi$ so that
\begin{eqnarray*}
&&c(x+1)^{-2}2\phi+(a_1+1)\phi_{xx}=0,\qquad\phantom{.aa}\phi(0)=a_2,\\
\noalign{\medskip}
&&c(x+1)^{-2}\phi^2-\textstyle\frac12\phi_x^2-(a_1+1)\psi_x=0,\quad\psi(0)=a_3\,.
\end{eqnarray*}
Set
$$T(x,y,\tilde x)=((a_1+1)x+a_1,y+\phi,(a_1+1)^{-1}\tilde x-\phi_xy+\psi)\,.$$
Then we may compute that:
\begin{eqnarray*}
&&T_*\partial_x=(a_1+1)\partial_x+\phi_x\partial_y+(-\phi_{xx}y+\psi_x)\partial_{\tilde x},\\
\noalign{\medskip}
&&T_*\partial_y=\partial_y-\phi_x\partial_{\tilde x},\qquad
T_*\partial_{\tilde x}=(a_1+1)^{-1}\partial_{\tilde x},
\end{eqnarray*}
and hence
\begin{eqnarray*}
&&g_f(T_*\partial_x,T_*\partial_x)=
-2\left\{(a_1+1)^2c((a_1+1)x+a_1+1)^{-2}(y^2+2\phi y+\phi^2)\right.\\
\noalign{\smallskip}
&&\phantom{g_f(T_*\partial_x,T_*\partial_x)=-2\left\{\right.}
\left.-\textstyle\frac12\phi_{x}^2
+(a_1+1)(\phi_{xx}y-\psi_x)\right\}=f,\\
\noalign{\medskip}
&& g_f(T_*\partial_x,T_*\partial_{\tilde x})=1,\qquad
g_f(T_*\partial_y,T_*\partial_y)=1,\\
\noalign{\medskip}
&&g_f(T_*\partial_x,T_*\partial_y)=0,\qquad
g_f(T_*\partial_y,T_*\partial_{\tilde x})=0,\qquad
g_f(T_*\partial_{\tilde x},T_*\partial_{\tilde x})=0\,.
\end{eqnarray*}
Thus $T$ is an isometry; since $T(0,0,0)=(a_1,a_2,a_3)$,
we have that $\mathcal{N}_c$
is a local homogeneous space. Furthermore, the shift
$T(x,y,\tilde x)=(x+a_1,y,\tilde x)$ provides an isometry between
$\mathcal{M}_{c(x+1)^{-2}}$ and $\mathcal{M}_{c(x+1+a_1)^{-2}}$
thereby showing these manifolds are isometric as well.
\end{remark}

\subsubsection{The model manifold $\mathcal{CW}_\varepsilon$ of Definition~\ref{def-2.3}}
Let $f(x,y)=\varepsilon y^2$ where $\varepsilon>0$.
Let $\beta=\beta(x)$ and $\gamma=\gamma(x)$ be given.
Choose $\phi$ so that $2\varepsilon\phi+\phi_{xx}=\beta$;
$\phi$ need not be globally defined, but this is always possible locally.
Then choose $\psi$ so
$-\psi_x-\frac12\phi_x^2+\varepsilon\phi^2=\gamma$.
Then
\begin{eqnarray*}
\tilde f(x,y)&=&\varepsilon y^2+(2\varepsilon\phi+\phi_{xx})y
-\psi_x-{\textstyle\frac12\phi_x^2}+\varepsilon\phi^2\\
\noalign{\medskip}
&=&\varepsilon y^2+\beta y+\gamma\,.
\end{eqnarray*}
Consequently, $T$ is a local isometry between $\mathcal{M}_{\varepsilon y^2}$
and $\mathcal{M}_{\varepsilon y^2+\beta y+\gamma}$. Since the transformation
$T_\varepsilon(x,y,\tilde x)=(\sqrt{\varepsilon} x,y,\tilde x/\sqrt\varepsilon)$
provides an isometry between $\mathcal{M}_{y^2}$ and
$\mathcal{M}_{\varepsilon y^2}$, the parameter $\varepsilon$ plays no role.
This completes the proof of Theorem~\ref{thm-1.4}~(1).

\begin{remark}\rm
We can use this formalism to see that $\mathcal{M}_{y^2}$ is a
homogeneous space (it is in fact a symmetric space). Suppose given a point
$(a_1,a_2,a_3)\in\mathbb{R}$.
Set $\phi(x)=a_2\cos(\sqrt{2}x)$. We then have that
$2\phi+\phi_{xx}=0$ and $\phi(0)=a_2$.
Now choose $\psi(x)$ so that $\psi_x+\frac12\phi_x^2+\phi^2=0$ and so that
$\psi(0)=a_3$. Let
$$T(x,y,z)=(x+a_1,y+\phi,\tilde x-\phi_xy+\psi)\,.$$
The translation in the $x$ coordinate
is harmless and does not change the equations of structure. We then have
that $T^*g_f=g_f$ and $T(0,0,0)=(a_1,a_2,a_3)$.
Consequently $\mathcal{M}_{y^2}$
is globally a homogeneous space.
\end{remark}

\begin{remark}\rm
Locally conformally flat homogeneous Lorentzian 3-manifolds with nilpotent Ricci operator were investigated in \cite{HT}. As any homogeneous Lorentzian three-manifold with nilpotent Ricci operator is a Walker manifold \cite{CLGRVAVL}, it follows from Theorem \ref{thm-1.4} that any non-symmetric  locally conformally flat homogeneous three-manifold with nilpotent Ricci operator is locally isometric to the manifold $\mathcal{P}_c$.

Further observe that all non-conformally flat left-invariant metrics on three-dimensional Lie groups
 (see, for example the discussion in \cite{Calv}, \cite{HaLee}) with nilpotent Ricci operators are locally isometric to the manifold $\mathcal{N}_b$.
\end{remark}

\begin{remark}\rm
There is a different notion of curvature homogeneity that is due
to Kowalski and Van\v{z}urov\'{a} \cite{KV}, \cite{KV2}.
Motivated by their seminal work, we say that a manifold $(M,g)$ is
\emph{Kowalski-Van\v{z}urov\'{a} $k$-curvature homogeneous} if for any two points there exists a linear
homothety between the corresponding tangent spaces which preserves the $(1,3)$-curvature {operator}
$\mathfrak{R}$ and its covariant derivatives up to order $k$. This concept lies between the notion of
{affine $k$-curvature} homogeneity and {$k$-curvature homogeneity}
since the {group of homotheties} lies between the orthogonal group and the general linear group.

It was shown in \cite{KV2} that the existence of a linear homothety
$\Phi_{p,q}:T_pM\rightarrow T_qM$ (i.e., $\Phi_{p,q}^*g_q=\lambda_{p,q}^2g_p$)
such that $\Phi_{p,q}^*\nabla^l\mathfrak{R}_q=\nabla^l\mathfrak{R}_p$ is equivalent to the existence of a
linear isometry $\varphi_{p,q}:T_pM\rightarrow T_qM$ and {$0\ne\lambda_{p,q}\in\mathbb{R}$}
{so}
${\lambda_{p,q}^{l+2}}\varphi_{p,q}^*\nabla^l{R}_q=\nabla^l{R}_p$.
Hence, $(M,g)$ is \emph{Kowalski-Van\v{z}urov\'{a} $k$-curvature homogeneous} if for any two points there exists a linear isometry
$\varphi_{p,q}$ between the corresponding tangent spaces such that
\begin{equation}\label{eq:KV1}
{\lambda_{p,q}}^{l+2}\varphi_{p,q}^*\nabla^l{R}_q
=\nabla^l{R}_p, \quad \mbox{for all}\quad l=0,\dots,k
\end{equation}
{or, equivalently, if there are constants $\varepsilon_{ij}$ and $c_{i_1...i_{\ell+4}}$ so that
for every point $p$ of $M$, there is a basis $\{\xi_1,...,\xi_m\}$ for $T_pM$ and a constant
$0\ne\lambda_p\in\mathbb{R}$ so that
\begin{equation}\label{eqn-4.b}
\begin{array}{l}
g(e_i,e_j)=\varepsilon_{ij},\text{ and}\\
\nabla^\ell R_p(e_{i_1},e_{i_2},e_{i_3},e_{i_4};e_{i_5}...e_{i_\ell+4})=\lambda_p^{2+\ell}
   c_{i_1...i_{\ell+4}}\text{ for all }\ell=0,...,k\,.
\vphantom{\vrule height 11pt}\end{array}\end{equation}
This is genuinely a different concept as the following illustrates.}

Let $\mathcal{M}_f$ be a Walker manifold with metric defined by  \eqref{eq:walker-metric}.
{We suppose that $f_{yy}>0$ and $f_{yyy}>0$.}
Proceed as in Section~\ref{section-3} by setting
$$
\xi_1:=a_{11}(\partial_x+f\partial_{\tilde x}+a_{12}\partial_y+a_{13}\partial_{\tilde x}),\quad
     \xi_2:=a_{22}\partial_y+a_{23}\partial_{\tilde x},\quad
     \xi_3:=a_{33}\partial_{\tilde x}\,.
$$
To ensure that $\{\xi_1,\xi_2,\xi_3\}$ is a pseudo-orthonormal {basis} (i.e., the non-zero components of the metric are given by $\langle \xi_1,\xi_3\rangle=\langle\xi_2,\xi_2\rangle=1$), we set:
$$
a_{22}=1, \quad a_{11}a_{33}=1,\quad a_{12}+a_{23}=0,
\quad a_{12}^2+2a_{13}=0\,.
$$
The crucial point is that $\{a_{11},a_{12}\}$ are free parameters and we use the relations above to determine
$a_{13}$, $a_{23}$, and $a_{33}$; these variables play no further role.
Now observe that
$$
\begin{array}{c}
\mathcal{R}(\xi_1,\xi_2,\xi_2,\xi_1)=a_{11}^2f_{yy},
\\
\noalign{\medskip}
\nabla\mathcal{R}(\xi_1,\xi_2,\xi_2,\xi_1;\xi_1)=a_{11}^3\left(f_{xyy}+a_{12}f_{yyy}\right),
\quad
\nabla\mathcal{R}(\xi_1,\xi_2,\xi_2,\xi_1;\xi_2)=a_{11}^2f_{yyy}\,.
\end{array}
$$
{We take $c_{1221}=1$, $c_{12211}=0$, and $c_{12212}=1$, in Equation~(\ref{eqn-4.b})
to define our model. This implies that
$$a_{11}^2f_{yy}=\lambda^2,\quad f_{xyy}+a_{12}f_{yyy}=0,\quad
\text{and}\quad a_{11}^2f_{yyy}=\lambda^3\,.$$
We solve these equations to obtain:
$$a_{12}=-f_{xyy}f_{yyy}^{-1},\quad
\lambda=f_{yyy}f_{yy}^{-1},\text{ and }a_{11}=\lambda f_{yy}^{-1/2}\,.$$
This shows} that \emph{any Walker metric \eqref{eq:walker-metric} such that $f_{yy}> 0$ and $f_{yyy}> 0$ is
Kowalski-Van\v{z}urov\'{a} 1-curvature homogeneous}. This is in contrast with the usual
{$1$-curvature homogeneity} (see Theorem \ref{thm-1.3}).

A more striking fact is the \emph{existence of three-manifolds which are Kowalski-Van\v{z}urov\'{a} curvature homogeneous of any order without being locally homogeneous}. In the Lorentzian case, such examples can be constructed by considering Walker manifolds $\mathcal{M}_f$, where $f(y)=\ln(y)$ (see \cite{GGN} for more details).\end{remark}

\section*{Acknowledgments}
Research of all the authors was partially supported by projects MTM2009-07756
and INCITE09 207 151 PR (Spain). Research of P. Gilkey and S. Nik\v cevi\'c was also
partially supported by project 174012 (Serbia).

%\bibliographystyle{model1a-num-names}



\end{document}